\documentclass[11pt,reqno]{amsart}
\usepackage{xifthen}
\usepackage{subfig}
\usepackage{pkgfile}

\newcommand{\bbV}{{\mathcal V}}
\newcommand{\bbI}{{\mathcal I}}
\newcommand{\bbH}{{\mathcal H}}

\newtheorem{obs}{Observation}
\newcommand{\remove}[1]{}

\begin{document}

\title[On Pseudo-Convex Partitions of a Planar Point Set]{On Pseudo-Convex Partitions of a Planar Point Set}

\author{Bhaswar B. Bhattacharya}
\address{Department of Statistics, Stanford University, California, USA, {\tt bhaswar@stanford.edu}}
   
\author{Sandip Das}
\address{Advanced Computing and Microelectronics Unit, Indian Statistical Institute, Kolkata, India, {\tt sandipdas@isical.ac.in}}

\begin{abstract}
Aichholzer et al. [{\it Graphs and Combinatorics}, Vol. 23, 481-507, 2007]
introduced the notion of pseudo-convex partitioning of planar point sets and
proved that the pseudo-convex partition number $\psi(n)$ satisfies,
$\frac{3}{4}\lfloor\frac{n}{4}\rfloor\leq \psi(n)\leq\lceil\frac{n}{4}\rceil$.
In this paper we prove that $\psi(13)=3$, which immediately improves the upper bound on $\psi(n)$ to $\lceil\frac{3n}{13}\rceil$, thus answering a question posed by Aichholzer et
al. in the same paper.
\end{abstract}


\maketitle

\section{Introduction}

In 1978 Erd\H os \cite{erdosempty} asked whether for every positive integer $k$,
there exists a smallest integer $H(k)$, such that any set of at least
$H(k)$ points in the plane, no three on a line, contains $k$ points which lie on the vertices of a convex polygon
whose interior contains no points of the set. Such a subset is called an {\it empty
convex $k$-gon} or a {\it k-hole}. Esther Klein showed {$ H(4)=5$} and Harborth \cite{harborth}
proved that {$ H(5)=10$}. Horton \cite{horton} showed
that it is possible to construct arbitrarily large set of points
without a 7-hole, proving that {$ H(k) $} does not exist
for {$ k \geq 7$}. Recently, after a long wait, the existence of $ H(6)$ has been proved
by Gerken \cite{gerken} and independently by Nicol\'as
\cite{nicolas}. Later Valtr \cite{valtrhexagon} gave a simpler version of Gerken's proof.

Any two empty convex polygons are said to be {\it disjoint} if their convex hulls
do not intersect. Let $H(k,\ell)$, $k\leq \ell$ denote the smallest integer
such that any set of $H(k, \ell)$ points in the plane, no three on a line,
contains both a $k$-hole and a $\ell$-hole which are disjoint. Clearly,
$H(3,3)=6$ and Horton's result \cite{horton} implies that $H(k, \ell)$ does not
exist for all $\ell\geq 7$.
It is known that $H(3, 4)=7$ \cite{urabedam}, $H(3,5)=10$ \cite{hosono},  $H(4,4)=9$  \cite{urabecgta}, and $H(4, 5)=12$ \cite{bbbsdqp}. Recently, Hosono and Urabe \cite{kyotocggt} showed that $H(5, 5) \geq 17$, and Bhattacharya and Das \cite{bbbsdpentagon} proved $H(5, 5) \leq 19$.

The problem of partitioning planar point sets with disjoint holes was first addressed by
Urabe \cite{urabedam}. For any set $S$ of points in the plane, denote the {\it convex hull} of $S$ by $CH(S)$, and cardinality of $S$ by $|S|$. Given a set $S$ of $n$ points in the plane, no three on a line, a {\it disjoint convex partition} of $S$ is a partition of $S$ into subsets $S_1, S_2, \ldots S_t$, with $\sum_{i=1}^t |S_i| = n$, such that for each $i\in\{1, 2, \ldots, t\}$, $CH(S_i)$ forms a $|S_i|$-gon and $CH(S_i)\cap CH(S_j)=\emptyset$, for any pair of distinct indices $i, j$. Observe that in any disjoint convex partition of $S$, the set $S_i$ forms a $|S_i|$-hole and the holes formed by the sets $S_i$ and $S_j$ are disjoint for any pair of distinct indices $i, j$. Let $\kappa(S)$ denote the minimum number of disjoint holes in any disjoint convex partition of $S$. Define
$\kappa(n)= \max_S\kappa(S)$, where the maximum is taken over all sets $S$ of $n$ points.
$\kappa(n)$ is called the {\it convex partition number} for all sets $S$ of fixed
size $n$, and it is bounded by $\lceil\frac{n-1}{4}\rceil \leq \kappa(n) \leq \lceil\frac{5n}{18}\rceil$.
The lower bound was given by Urabe \cite{urabedam} and the upper bound
by Hosono and Urabe \cite{urabecgta}. The lower bound was later improved to $\lceil\frac{n+1}{4}\rceil$ by Xu and Ding \cite{empty_c_annals}.


A pseudo-triangle is a simple polygon with exactly three vertices having interior angles less than
$180^{\circ}$, and is considered to be the natural counterpart of a convex polygon. These have been studied recently in the context of pseudo-triangulations, which are tessellation of the plane with pseudo-triangles. They provide sparser tessellations than triangulations but retain many of the desirable properties of triangulations. Pseudo-triangulations have received considerable attention in the last few years for applications in areas like motion planning, collision detection, ray shooting, rigidity, or visibility (refer to Rote et al. \cite{surveypseudo} for a survey of the different properties of pseudo-triangulations and its various applications).

A pseudo-triangle with $\ell$ vertices is called a $\ell$-pseudo-triangle, and a set is said to contain an empty $\ell$-pseudo-triangle if there exists a subset of $\ell$ points forming a pseudo-triangle which contains
no point of the set in its interior. Any two empty pseudo-triangles, or a hole and an empty pseudo-triangle are said to be disjoint if their vertex sets as well as their interiors are disjoint. Recently, Aichholzer et al. \cite{toth} introduced the problem of partitioning planar point sets with disjoint holes or empty pseudo-triangles. Given a set $S$ of $n$ points in the plane, no three on a line, a {\it pseudo-convex partition} of $S$ is a partition of $S$ into subsets $S_1, S_2, \ldots S_t$, with $\sum_{i=1}^t |S_i| = n$, such that for each $i\in\{1, 2, \ldots, t\}$, the set $S_i$ forms a $|S_i|$-hole or a $|S_i|$-pseudo-triangle, the holes or pseudo-triangles formed by the sets $S_i$ and $S_j$ are disjoint for any pair of distinct indices $i, j$. If $\psi(S)$ denotes the minimum number of disjoint holes or empty pseudo-triangles in any pseudo convex partition of $S$, then the {\it pseudo-convex partition number} is defined as $\psi(n)= \max_S\psi(S)$, where the maximum is taken over all sets $S$ of $n$ points.

Aichholzer et al. \cite{toth} showed that the pseudo-convex partition number $\psi(n)$
satisfies: $\frac{3}{4}\lfloor\frac{n}{4}\rfloor\leq \psi(n)\leq\lceil\frac{n}{4}\rceil$.
The upper bound follows from the simple observation that every set of 4 points forms either an
4-hole or an empty 4-pseudo-triangle. In fact, using computer-aided search, Aichholzer et al. \cite{toth} obtained bounds on the pseudo-convex partition number $\psi(n)$ for small point sets. However, they were unable to find the exact value of $\psi(13)$, and mentioned the possibility of a non-trivial upper bound on $\psi(n)$ by conjecturing that $\psi(13)=3$. In this paper, we answer this question in the affirmative, thus proving that $\psi(n)\leq\lceil\frac{3n}{13}\rceil$. Our proof is geometric and does not rely on computer-aided search over the order type database. We identify a number of simple necessary conditions that allows the desired partitioning, and then proceed to show that each set of 13 points must satisfies one of these conditions.

\section{Notations and Definitions}
\label{sec:notations}

We first introduce the definitions and notations required for the remaining part of the paper.
Let $S$ be a finite set of points in the plane in general position, that is, no three on a line.
Denote the {\it convex hull} of $S$ by $CH(S)$. The
boundary vertices of $CH(S)$, and the points of $S$ in the interior
of $CH(S)$ are denoted by $\bbV(CH(S))$ and $\tilde{\bbI}(CH(S))$, respectively. A region $R$ in the plane is said to be {\it empty} in $S$ if $R$ contains no elements of $S$ in its interior. Moreover, for any set $T$, $|T|$ denotes the cardinality of $T$.

By $\mathcal P:=p_1p_2\ldots p_m$ we denote the region bounded by the simple polygon with vertices $\{p_1, p_2, \ldots, p_m\}$ ordered anti-clockwise. Let $\bbV(\mathcal P)$ denote the set of vertices $\{p_1, p_2, \ldots, p_m\}$ and $\bbI(\mathcal P)$ the interior of $\mathcal P$. A finite set of points $Z$ is said to {\it span} a simple polygon $\mathcal P$ if $\bbV(\mathcal P)=Z$.

\begin{figure*}[h]
\centering
\begin{minipage}[l]{0.5\textwidth}
\centering
\includegraphics[width=2.95in]
    {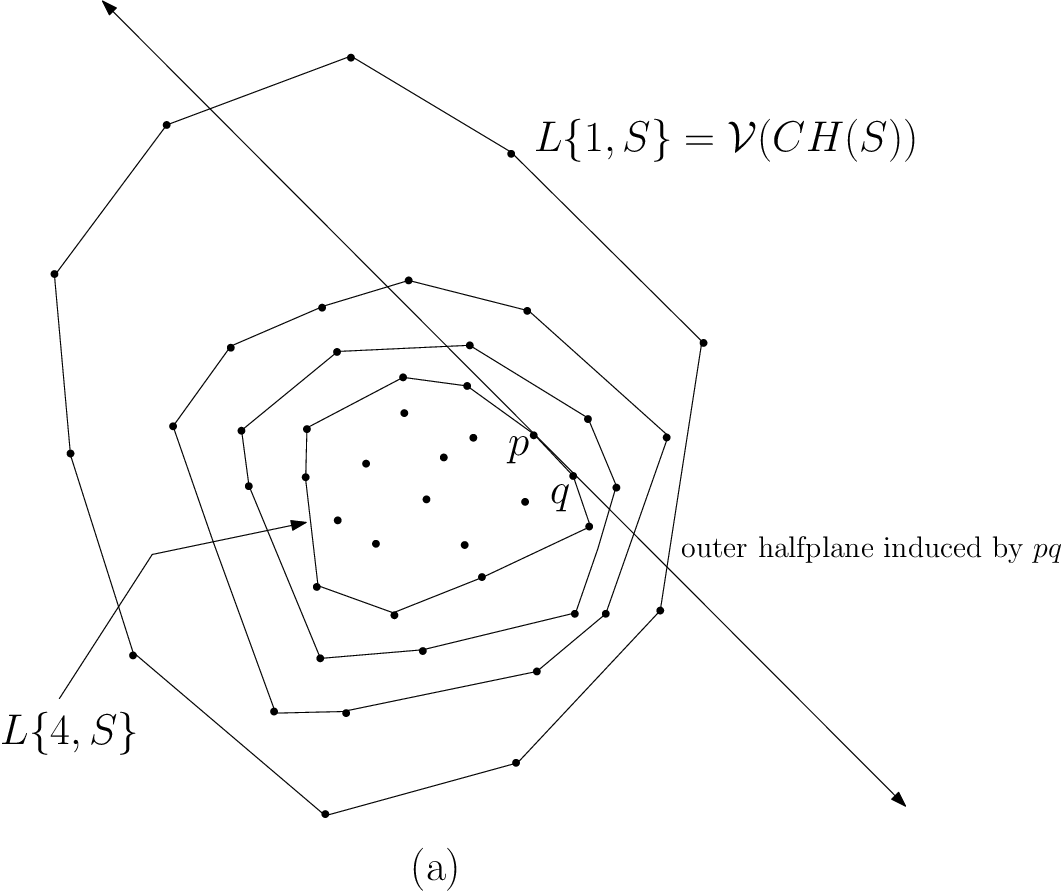}\\
\end{minipage}
\begin{minipage}[c]{0.4\textwidth}
\centering
\includegraphics[width=2.95in]
    {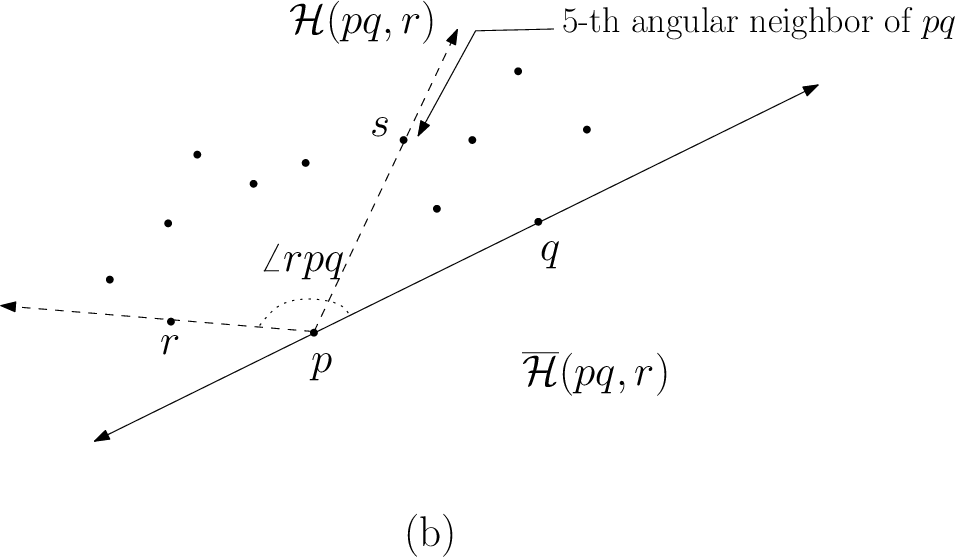}\\
\end{minipage}%
\caption{(a) Convex layers and outer halfplane, and (b) Cone and nearest angular neighbors.}
\label{fig:definitions}
\end{figure*}

The {\it $j$-th convex layer} of $S$, denoted by $L\{j, S\}$, is the
set of points  that lie on the boundary of
$CH(S\backslash\{\bigcup_{i=1}^{j-1}L\{i, S\}\})$, where $L\{1, S\}=\bbV(CH(S))$.
If $p,~q\in S$ are such that ${pq}$ is an edge of the convex
hull of the $j$-th layer, then the open halfplane bounded by the line
$pq$ and not containing any point of $S\backslash\{\bigcup_{i=1}^{j-1}L\{i, S\}\}$ will be referred to as the {\it outer} halfplane induced by the edge ${pq}$ (see Figure \ref{fig:definitions}(a)).

For any three points $p, q, r \in S$, $\bbH(pq, r)$ (respectively $\mathcal H_c(pq, r)$) denotes the open (respectively closed) halfplane bounded by the line $pq$ containing the point $r$. Similarly, $\overline{\bbH}(pq, r)$ (respectively $\overline{\bbH}_c(pq, r)$) is the open (respectively closed) halfplane bounded by $pq$ not containing the point $r$. Moreover, if $\angle rpq < \pi$, $Cone(rpq)$ denotes the interior of the angular domain $\angle rpq$. A point $s \in Cone(rpq)\cap S$ is called the {\it nearest angular neighbor} of $\overrightarrow{pq}$ in $Cone(rpq)$ if $Cone(spq)$ is empty in $S$. Similarly, for any convex region $R$ a point $s \in R\cap S$ is called the {\it nearest angular neighbor} of $\overrightarrow{pq}$ in $R$ if $Cone(spq)\cap R$ is empty in $S$. More generally, for any positive integer
$k$, a point $s\in S$ is called the {\it $k$-th angular neighbor} of $\overrightarrow{pq}$ whenever $Cone(spq)\cap R$ contains exactly $k-1$ points of $S$ in its interior (see Figure \ref{fig:definitions}(b)).

\section{Pseudo-Convex Partitioning}
\label{sec:theorem_proof}

Aichholzer et al. \cite{toth} showed that $\psi(n) \leq \lceil\frac{n}{4}\rceil$. They also observed that $3\leq \psi(13)\leq 4$, and mention the possibility of a better upper bound of $\lceil\frac{3n}{13}\rceil$ on $\psi(n)$
by conjecturing that $\psi(13)=3$.

In the following theorem we settle this conjecture in the affirmative.

\begin{thm}
Every set of 13 points in the plane, in general position, can be partitioned into
three sets each of which span either a hole or
an empty pseudo-triangle which are mutually disjoint. In other words, $\psi(13)=3$.
\label{lm3}
\end{thm}

Theorem \ref{lm3} immediately establishes a non-trivial upper bound on $\psi(n)$, as suggested by Aichholzer et al. \cite{toth}:

\begin{thm}
$\psi(n)\leq \lceil\frac{3n}{13}\rceil.$
\label{th5}
\end{thm}

\begin{proof}Let $S$ be a set of $n$ points in the plane, no three of which are collinear.
By a horizontal sweep we can divide the plane into
$\lceil\frac{n}{13}\rceil$ disjoint strips, of which $\lfloor\frac{n}{13}\rfloor$ contain 13 points each
and one remaining strip $R$, with $|R|<13$.
The strips having 13 points, can be partitioned into
three disjoint holes or empty pseudo-triangles by Theorem \ref{lm3}.
Since $|R|<13$, at most $\lceil\frac{3|R|}{13}\rceil$ disjoint holes or empty pseudo-triangles are needed to partition $R$,
thus proving that $\psi(n)\leq \lceil\frac{3n}{13}\rceil$. 
\end{proof}


\section{Proof of Theorem \ref{lm3}}

Let $S$ be a set of 13 points in the plane in general position. A partition of $S$ into three disjoint subsets $S_1, S_2, S_3$ is called {\it admissible} if each $S_i$, $i\in\{1, 2, 3\}$, is either empty or it forms an $|S_i|$-hole or empty $|S_i|$-pseudo-triangle, such that the holes or pseudo-triangles formed by the sets $S_i$ and $S_j$ are disjoint for any pair of distinct indices $i\ne j$.
The set $S$ is said to be {\it admissible} if there exists an admissible partition of $S$.
To prove Theorem \ref{lm3} we need to exhibit an admissible partition of $S$, for all sets of 13 points in the plane, in general position.

Observe that any set of 4 points in the plane always spans a convex quadrilateral or a 4-pseudo-triangle.
Hence, we have the following observation.

\begin{obs}For every integer $k\geq 1$, we have $\psi(4k)\leq k$. \hfill $\Box$
\label{ob:psi4k}
\end{obs}

\begin{obs}
$S$ is admissible if some outer halfplane induced by an edge of the second convex layer
contains more than two points of $\bbV(CH(S))$.
\label{ob:ob1}
\end{obs}

\begin{proof}Suppose some outer halfplane induced by an edge of the second layer
contains more than two points of $\bbV(CH(S))$. This means that there exists two points $p, q\in S$, such
that the line segment $pq$ is an edge of the second convex layer and $|\overline{\bbH}(pq, r)\cap S| \geq 3$,
where $r \in L\{2, S\}\backslash \{p, q\}$. Then $\overline{\bbH}_c(pq, r)\cap S$ spans a $k$-hole, with $k \geq 5$. The remaining eight points of $S$ all lie in the halfplane $\mathcal H(pq, r)$. As $|\mathcal H(pq, r)\cap S|\leq 8$, the points in $\mathcal H(pq, r)\cap S$ can be partitioned using at most two disjoint holes or empty pseudo-triangles. 
\end{proof}

\begin{figure*}[h]
\centering
\begin{minipage}[c]{0.33\textwidth}
\centering
\includegraphics[width=1.90in]
    {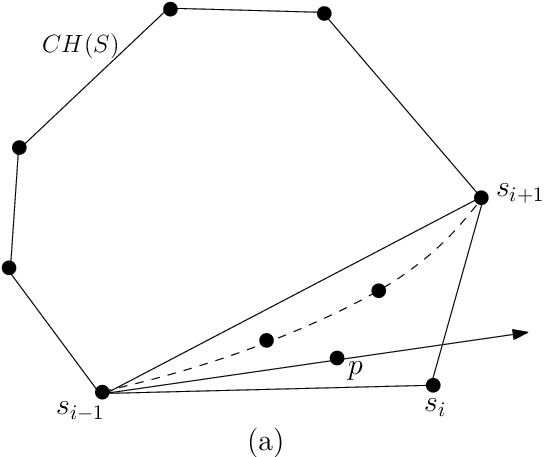}\\
\end{minipage}%
\begin{minipage}[c]{0.33\textwidth}
\centering
\includegraphics[width=1.90in]
    {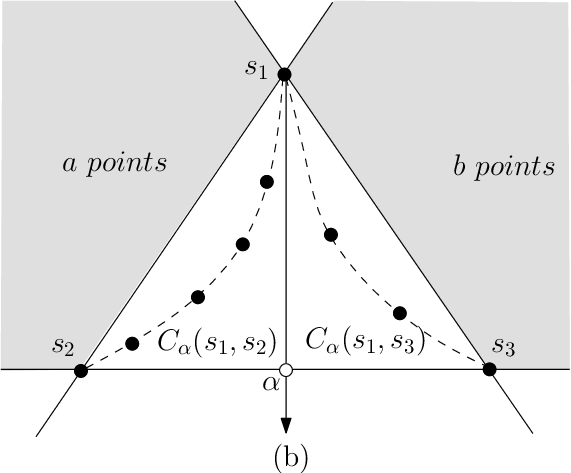}\\
\end{minipage}
\begin{minipage}[c]{0.33\textwidth}
\centering
\includegraphics[width=1.90in]
    {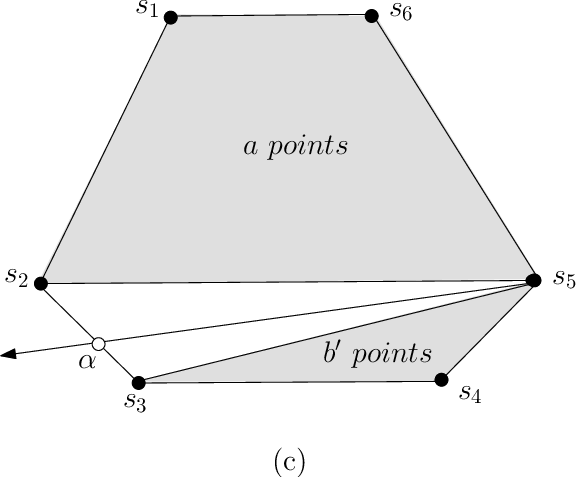}\\
\end{minipage}%
\caption{Illustration for the proof of (a) Observation \ref{ob:ob2}, (b) Observation \ref{ob:obnew}, and (c) Observation \ref{obs:splitter}.}
\label{figure1}
\end{figure*}

\begin{obs}Let $\bbV(CH(S))=\{s_1, s_2,\dots, s_k\}$, with the vertices in counter-clockwise order. If there exists a point $s_i\in \bbV(CH(S))$, such that $|\bbI(s_{i-1}s_{i}s_{i+1})\cap S|\geq 2$, then $S$ is admissible,
where the indices are taken modulo $k$.
\label{ob:ob2}
\end{obs}

\begin{proof}Suppose $|\bbI( s_{i-1}s_is_{i+1})\cap S|=a\geq 2$, for some $i=1,2, \ldots ,k$.
Let $p\in S$ be the first angular neighbor of $\overrightarrow{s_{i-1}s_i}$ in $Cone(s_is_{i-1}s_{i+1})$.
Let $Z=(\mathcal H_c(s_{i-1}s_{i+1}, s_i)\cap S)\backslash \{p, s_i\}$, and suppose $\bbV(CH(Z))=\{s_{i-1}, z_1, \ldots, z_b, s_{i+1}\}$. As $a\geq 2$, $b \geq 1$, and $s_{i+1}z_1\ldots z_b s_{i-1}ps_i$ is an empty $\ell$-pseudo-triangle with $\ell \geq 5$. Let $S'=S\backslash \{s_{i+1}, z_1, \ldots, z_b, s_{i-1}, p, s_i\}$, that is, the points of $S$ without the points of the considered $\ell$-pseudo-triangle. As $|S'|\leq 8$ and $CH(S')$ does not intersect the $\ell$-pseudo-triangle, $S'$ can be partitioned using at most two holes or empty pseudo-triangles from Observation \ref{ob:psi4k}. 
\end{proof}

\begin{obs}$S$ is admissible if there exists three distinct points $s_1,s_2,s_3\in S$ satisfying the following two conditions:
\item[(A)]$|(\overline{\bbH}(s_1s_2, s_3)\cap \overline{\bbH}(s_1s_3, s_2))\cap S|=0$ and $|\overline{\bbH}(s_2s_3, s_1)\cap S|=0$.
\item[(B)]{$|\overline{\bbH}(s_1s_2, s_3)\cap S|\leq 4$ and $|\overline{\bbH}(s_1s_3, s_2)\cap S|\leq 4$.}
\label{ob:obnew}
\end{obs}

\begin{proof}Let $|\overline{\bbH}(s_1s_2, s_3)\cap S|=a\leq 4$ and $|\overline{\bbH}(s_1s_3, s_2)\cap S|=b\leq 4$.
Then, $|\bbI(s_1s_2s_3)\cap S|=10-(a+b)$. Then, there exists a point $\alpha \notin S$ on the line segment ${s_2s_3}$ such that
$|\bbI( s_1s_2\alpha)\cap S|=5-a$ and $|\bbI( s_1s_3\alpha)\cap S|=5-b$ (see Figure \ref{figure1}(b)).
Let $C_\alpha(s_1,s_2)=CH(\{\bbI(s_1s_2\alpha)\cap S\}\cup\{s_1, s_2\})$ and
$C_\alpha(s_1, s_3)=CH(\{\bbI(s_1s_3\alpha)\cap S\}\cup\{s_1, s_3\})$.
Now, since both $a, b\leq 4$, we have $|\bbV(C_\alpha(s_1, s_2))|\geq 3$ and $|\bbV(C_\alpha(s_1, s_3))|\geq 3$. Thus,
$S_1=\bbV(C_\alpha(s_1, s_2))\cup\bbV(C_\alpha(s_1, s_3))$ spans an empty $\ell$-pseudo-triangle, with $\ell\geq 5$.
Moreover, both $S_2=\bbI(C_\alpha(s_1, s_2))\cup \{\overline{\bbH}(s_1s_2, s_3)\cap S\}$, and $S_3=\bbI(C_\alpha(s_1, s_3))\cup\{\overline{\bbH}(s_1s_3, s_2)\cap S\}$ lie inside two disjoint convex regions containing at most 4 points each. Therefore, the partition $S=S_1\cup S_2\cup S_3$ is admissible. 
\end{proof}

If there exist three distinct points $s_1, s_2, s_3\in S$ satisfying conditions (A)
and (B) of Observation \ref{ob:obnew}, then the three points are said to
form a {\it heart} with the line segment $s_2s_3$ as {\it base} and the point $s_1$ as {\it pivot}.
The set $S$ is admissible if any three of its points form a {\it heart}. Observe that if $s_1, s_2, s_3\in \bbV(CH(S))$ and some edge of the triangle $ s_1s_2s_3$ is also an edge of
$CH(S)$, then condition (A) is automatically satisfied. In such cases, $ s_1, s_2, s_3$ form a {\it heart} whenever condition (B) holds.

Equipped with the above three observations, we now proceed to prove the admissibility of $S$.
The proof of the admissibility of $S$ is presented in three separate sections. The first section deals with the cases $|CH(S)|\leq 5$, the second section with the case $|CH(S)|=6$, and the third considers the cases $|CH(S)|\geq 7$.

Let $\bbV(CH(S))=\{s_1, s_2,\dots, s_k\}$, with the vertices taken in the counter-clockwise order.
While indexing a set of points from $\bbV(CH(S))$, we identify indices modulo $k$.

\subsection{$|CH(S)|\leq 5$}

\begin{lem}$S$ is admissible whenever $|CH(S)|\leq 5$.
\label{lm:chleq5}
\end{lem}

\begin{proof}Observe that if $|CH(S)|=3$, then the admissibility of $S$ is a direct consequence of
Observation \ref{ob:ob2}. Now, we consider the following two cases based on the size of $|CH(S)|$:

\begin{description}
\item[{\it Case} 1] $|CH(S)|=4$. This implies that $|\tilde{\bbI}(CH(S)|=9$. Therefore,
$|\bbI(s_2s_3s_4)\cap S|\geq 2$ or $|\bbI(s_1s_2s_4)\cap S|\geq 2$. The admissibility of $S$ then follows from Observation \ref{ob:ob2}.

\item[{\it Case} 2] $|CH(S)|=5$. Suppose that $|\bbI( s_1s_2s_3)\cap S|=a$, and
$|\bbI( s_1s_4s_5)\cap S|=b$. If $a \geq 2$ or $b \geq 2$, the admissibility of
$S$ is guaranteed from Observation \ref{ob:ob2}.
Therefore, assume that both $a, b\leq 1$. This implies that
$|\overline{\bbH}(s_1s_3, s_4)\cap S|\leq 2$ and $|\overline{\bbH}(s_1s_4, s_3)\cap S|\leq 2$. Thus, the three points $s_1, s_3, s_4$ satisfy Conditions (A) and (B) of Observation \ref{ob:obnew} and form a {\it heart} with ${s_3s_4}$ as {\it base} and $s_1$ as {\it pivot}. Thus, the admissibility of $S$ follows.
\end{description}
\end{proof}

\subsection{$|CH(S)|=6$}

For any $i\in\{1, 2, 3\}$, the diagonal $d:=s_is_{i+3}$ of the hexagon $s_1s_2s_3s_4s_5s_6$ is called an $(a, b)-splitter$ of $CH(S)$, where $a\leq b$ are integers, if either $|\mathcal H(s_is_{i+3}, s_{i+1}) \cap \tilde{\bbI}(CH(S))|=a$ and $|\overline{\bbH}(s_is_{i+3}, s_{i+1})\cap \tilde{\bbI}(CH(S))|=b$, or $|\overline{\bbH}(s_is_{i+3}, s_{i+1}) \cap \tilde{\bbI}(CH(S))|=a$ and $|\mathcal H(s_is_{i+3}, s_{i+1})\cap \tilde{\bbI}(CH(S))|=b$.

We now have the following observation:

\begin{obs}
If any one of the three diagonal ${s_2s_5}$, ${s_1s_4}$, and ${s_3s_6}$ is not a $(3, 4)$-splitter of $CH(S)$, then $S$ is admissible.
\label{obs:splitter}
\end{obs}

\begin{proof}It suffices to prove that $S$ is admissible whenever ${s_2s_5}$ is not a
$(3, 4)$-splitter of $CH(S)$. Suppose the diagonal ${s_2s_5}$ is a $(a, b)$-splitter of $CH(S)$, with
$a\leq 2$, and $b\geq 5$. Refer to Figure \ref{figure1}(c). W.l.o.g. assume that $|\bbI(s_1s_2s_5s_6)\cap S|=a$
and $|\bbI( s_2s_3s_4s_5)\cap S|=b$. Now, since $b\geq 5$, by the pigeon-hole principle, either $\bbI(s_2s_3s_5)$ or
$\bbI(s_3s_4s_5)$ contains at least 3 points of $S$. However, if $|\bbI( s_3s_4s_5)\cap S|\geq 2$, then
Observation \ref{ob:ob2} guarantees the admissibility of $S$.
Therefore, assume that $|\bbI(s_3s_4s_5)\cap S|=b'\leq 1$.
This implies that $|\overline{\bbH}(s_5s_3, s_2)\cap S|\leq 2$ and
$|\overline{\bbH}(s_5s_2, s_3)\cap S|\leq 4$, and the three points $s_2, s_3, s_5$ form a {\it heart} with $s_2s_3$ as {\it base} and $s_5$ as {\it pivot} (see Figure \ref{figure1}(c)). The admissibility of $S$ thus follows from Observation \ref{ob:obnew}. 
\end{proof}

In light of Observation \ref{obs:splitter}, it suffices to assume that the three diagonals ${s_2s_5}$, ${s_1s_4}$, and ${s_3s_6}$
are $(3, 4)$-splitters of $CH(S)$. Consider the partition of the interior of $CH(S)$ by the three diagonals into 7 disjoint
regions $R_i$ as shown in Figure \ref{fig6}(a). Let $|R_i|$ denote
the number of points of $S$ inside region $R_i$.

Now, we have the following observation:
\begin{figure*}[h]
\centering
\begin{minipage}[c]{0.5\textwidth}
\centering
\includegraphics[width=2.25in]
    {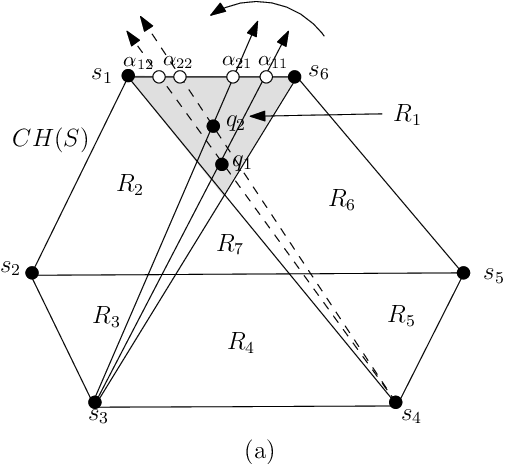}\\
\end{minipage}%
\begin{minipage}[c]{0.5\textwidth}
\centering
\includegraphics[width=2.25in]
    {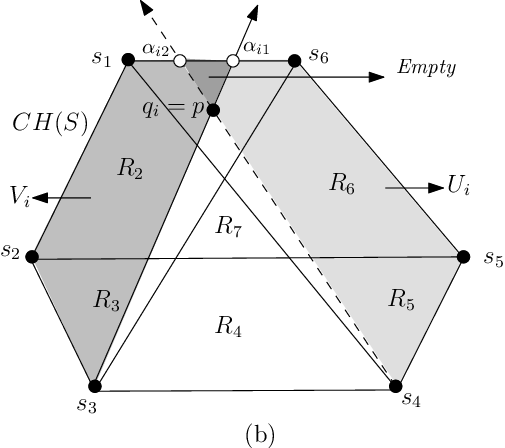}\\
\end{minipage}
\caption{Every diagonal of $CH(S)$ forms a $(3,4)$-splitter of the hexagon $s_1s_2s_3s_4s_5s_6$: (a) Illustration for the proof of Observation \ref{ob3}, (b) The points $ p, s_3, s_4$ form a {\it heart} with ${s_3s_4}$ as {\it base} and $p$ as pivot.}
\label{fig6}
\end{figure*}

\begin{obs}If $|R_4|+|R_7|\geq 2$ and $|R_1|\geq 1$, then there exists a point $p\in R_1\cap S$ such that
the three points $p, s_3, s_4$ form a heart with $s_3s_4$ as base with $p$ as pivot.
\label{ob3}
\end{obs}

\begin{proof}Let $|R_4|+|R_7|=a\geq 2$. Now, since both the diagonals $s_1s_4$ and $s_3s_6$ are $(3, 4)$-splitters of $CH(S)$,
we have $|R_2|+|R_3|=b\leq 2$, $|R_5|+|R_6|=b'\leq 2$, and $|R_1|+|R_5|+|R_6|\geq 3$. Let $q_1\in R_1\cap S$ be the $(3-b)$-th angular neighbor of $\overrightarrow{s_3s_6}$ in $R_1$. Let $U_1=(Cone(q_1s_3s_4)\backslash\bbI(q_1s_3s_4))\cap S$ and $V_1=\mathcal H(s_3q_1, s_2)\cap S$, and $\alpha_{11}$ and $\alpha_{12}$ are the points where the rays $\overrightarrow{s_3q_1}$ and
$\overrightarrow{s_4q_1}$ intersect the boundary $CH(S)$, respectively (Figure \ref{fig6}(a)). Therefore, $|U_1|\leq 4$ and $|V_1|\leq 4$. Now, if $Cone(\alpha_{11}q_1\alpha_{12})\cap S$ is empty, then $q_1(=p)$ and the result follows.

Otherwise, suppose that $Cone(\alpha_{11}q_1\alpha_{12})\cap S$ is non-empty.
Let $q_2\in R_1\cap S$ be the nearest angular neighbor of $\overrightarrow{s_3q_1}$ in $Cone(\alpha_{11}q_1\alpha_{12})$. Let $\alpha_{21}$ and $\alpha_{22}$
be the points where the rays $\overrightarrow{s_3q_2}$ and
$\overrightarrow{s_4q_2}$ intersect the boundary $CH(S)$, respectively.
Define, $U_2=(Cone(q_2s_3s_4)\backslash\bbI(q_2s_3s_4))\cap S$ and $V_2=\mathcal H(s_3q_2, s_2)\cap S$.
Observe that $U_2\subseteq U_1$ and $V_2\subseteq V_1$, and hence, $|U_2|\leq 4$ and $|V_2|\leq 4$.
Therefore, if $Cone(\alpha_{21}q_2\alpha_{22})\cap S$ is empty, then $q_2(=p)$ is the required point.

If $Cone(\alpha_{21}q_2\alpha_{22})\cap S$ is non-empty, we repeat the same procedure again, until we get a point $p~(=q_i)\in R_1\cap S$ with $|U_i|\leq 4$, $|V_i|\leq 4$, and $|Cone(\alpha_{i1}q_i\alpha_{i2})\cap S|=0$, where $q_i$ is the nearest angular neighbor of $\overrightarrow{s_3q_{i-1}}$ in $R_1\cap S$, $U_i=(Cone(q_is_3s_4)\backslash\bbI(q_is_3s_4))\cap S$, $V_i=\mathcal H(s_3q_i, s_2)\cap S$, and $\alpha_{i1}$ and $\alpha_{i2}$ are the points where the rays $\overrightarrow{s_3q_i}$ and $\overrightarrow{s_4q_i}$ intersect the boundary $CH(S)$, respectively (see Figure \ref{fig6}(b)). 
\end{proof}

The admissibility of $S$ when $|CH(S)|=6$ is now proved by considering the following three cases:

\begin{description}
\item[$Case$ 1] $|R_4|+|R_7|\geq 2$ and $|R_1|\geq 1$. In this case,
Observation \ref{ob3} guarantees the existence of a point $p\in
R_1\in S$ such that the three points $p, s_3, s_4$ form a
{\it heart} with ${s_3s_4}$ as {\it base} and $p$ as {\it pivot} (see Figure \ref{fig6}(b)).

\item[$Case$ 2] $|R_4|+|R_7|\geq 2$ and $|R_1|=0$. This implies that $|R_2|+|R_3|\leq 2$,
since ${s_1s_4}$ is a $(3,4)$-splitter of $CH(S)$.
Thus, $|R_1|+|R_2|+|R_3|\leq 2$, which contradicts the assumption that diagonal ${s_3s_6}$ is
a $(3, 4)$-splitter of $CH(S)$.

\item[$Case$ 3] If the previous two cases do not hold, then by symmetry we must have $|R_2|+|R_7|\leq 1$, $|R_4|+|R_7|\leq 1$, and $|R_6|+|R_7|\leq 1$. Therefore,
$|R_2|+|R_4|+|R_6|+|R_7|\leq 3$ which implies that $|R_1|+|R_3|+|R_5|\geq 4$.
By the pigeon-hole principle, one of the three regions $R_1$, $R_3$, and $R_5$ contains at least
two points $S$. W.l.o.g., assume $|R_1|\geq 2$. This implies that $|R_2|+|R_3|\leq 2$, and
hence $|R_4|+|R_7|\geq 1$. Combining this with the given inequality we get, $|R_4|+|R_7|=1$, $|R_2|+|R_3|=2$, and $|R_1|=2$. For $\{p\}\in (R_4\cup R_7)\cap S$, the three points $s_1, p, s_6$ form a {\it heart} with ${s_1s_6}$ as {\it base} and $p$ as {\it pivot}.
\end{description}

\subsection{$|CH(S)|=7$}

In this section we proof the admissibility of $S$ when $|CH(S)|=7$. As before assume that $L\{2, S\}=\{p_1, p_2, \ldots, p_m\}$, where the vertices are taken in counter-clockwise order, and the indices are to be identified modulo $m$.

\begin{lem}$S$ is admissible whenever $|CH(S)|=7$.
\label{lm:lmch7}
\end{lem}

\begin{proof}If $|L\{2, S\}|=3$, the admissibility of $S$ follows easily from Observation \ref{ob:ob1}.
Therefore, $4\leq |L\{2, S\}|\leq 6$. We consider these three cases separately as follows:

\begin{figure*}[h]
\centering
\begin{minipage}[c]{0.33\textwidth}
\centering
\includegraphics[width=2.15in]
    {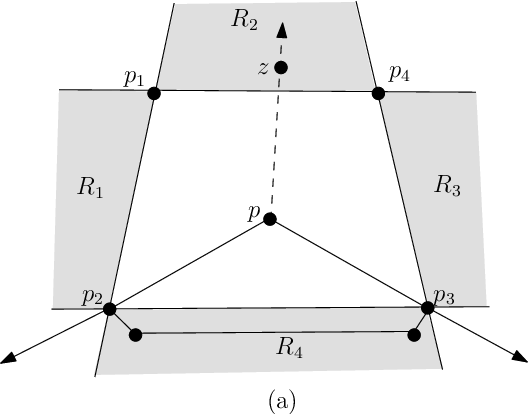}\\
\end{minipage}%
\begin{minipage}[c]{0.33\textwidth}
\centering
\includegraphics[width=2.0in]
    {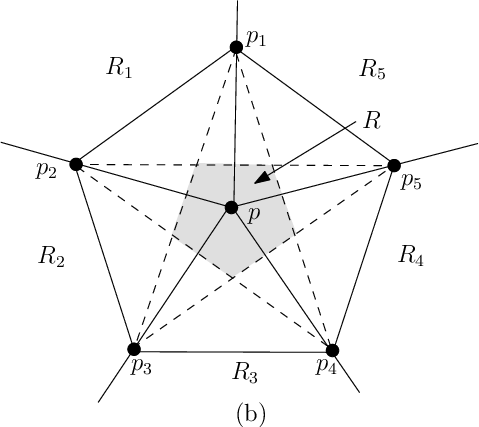}\\
\end{minipage}%
\begin{minipage}[c]{0.33\textwidth}
\centering
\includegraphics[width=1.95in]
    {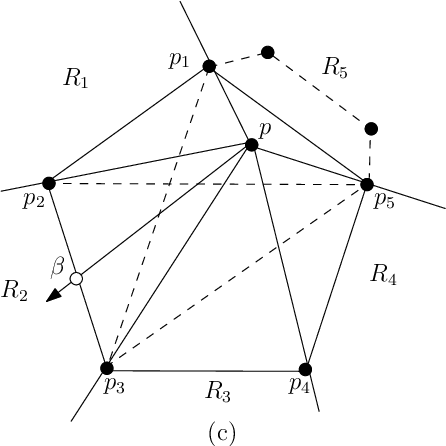}\\
\end{minipage}
\caption{Illustration for the proof of Lemma \ref{lm:lmch7}: (a) $|L\{2, S\}|=4$, (b) $|L\{2, S\}|=5$ and $p\in R$, and(c) $|L\{2, S\}|=5$ and $p\notin R$.}
\label{fig8}
\end{figure*}

\begin{description}
\item[{\it Case} 1]$|L\{2, S\}|=4$. Let $R_1, R_2, R_3, R_4$ be the 4 shaded regions outside the second convex layer, as shown in Figure \ref{fig8}(a). Note that if $|R_1|+|R_3|\leq 2$, then $S_2=((R_1\cup R_3)\cap S)\cup L\{3, S\}$, and forms a convex quadrilateral or a 4-pseudo-triangle. Hence, $S_1=\overline{\bbH}_c(p_1p_4, p_2)\cap S$, $S_2$, and $S_3=\overline{\bbH}_c(p_2p_3, p_1)\cap S$ is an admissible partition of $S$. Therefore, we may assume that $|R_1|+|R_3|\geq 3$, and, by symmetry, $|R_2|+|R_4|\geq 3$. As $\sum_{i=1}^{4}|R_i|\leq 7$, one of the above inequalities must be an equality. W.l.o.g. assume $|R_2|+|R_4|=3$. Moreover, by Observation \ref{ob:ob1} it suffices to assume that $|R_2|=1$ and $|R_4|=2$, and let $p\in L\{3, S\}$, be the nearest angular neighbor of $\overrightarrow{p_2p_3}$ inside the second convex layer. Then $S_2=(R_4\cap S)\cup \{p_2, p_3, p\}$ spans a 5-hole. If $R_2\cap S=\{z\}$, then depending on the position of the point $q\in L\{3, S\}$, either $|Cone(zpp_2)\cap S|=4$ or $|Cone(zpp_3)\cap S|=4$. Therefore, the partition $S_1=Cone(zpp_2)\cap S$, $S_2$ and $S_3=(Cone(zpp_3)\cap S)\cup\{z\}$ or the partition $S_1=(Cone(zpp_2)\cap S)\cup\{z\}$, $S_2$ and $S_3=(Cone(zpp_3)\cap S)$ is admissible for $S$, respectively.

\item[{\it Case} 2] $|L\{2, S\}|=5$. Let $L\{3, S\}=\{p\}$ and consider the partition of the exterior of the second convex layer into disjoint regions $R_i=Cone(p_ipp_{i+1})\backslash \bbI(p_ipp_{i+1})$, for $i\in\{1, 2, 3, 4, 5\}$. Let $R$ be shaded region inside the second convex layer as shown in Figure \ref{fig8}(b).

\begin{description}
\item[{\it Case} 2.1] $p\in R$. Observe that $\sum_{i=1}^5|R_i|=7$ and by Observation \ref{ob:ob1} for every $i\in \{1, 2, 3, 4, 5\}$, $|R_i|\leq 2$. Therefore, w.l.o.g. assume $|R_1|=2$ (Figure \ref{fig8}(b)).
If $|R_2|+|R_3|=2$ and $|R_4|+|R_5|=3$, then $S_1=(R_1\cap S)\cup\{p, p_1, p_2\}$, $S_2=Cone(p_1pp_4)\cap S$, and $S_3=(Cone(p_2pp_4)\cap S)\cup\{p_4\}$ is an admissible partition of $S$. Similarly, for $|R_2|+|R_3|=3$ and $|R_4|+|R_5|=2$. Otherwise, either $R_2\cup R_3$ or $R_4\cup R_5$ has more than 3 points of $S$. W.l.o.g., assume that $|R_2|+|R_3|\geq 4$. Observation \ref{ob:ob1} implies that $|R_2|=|R_3|=2$, and either the partition $S_1=(R_2\cap S)\cup\{p, p_2, p_3\}$, $S_2=Cone((p_2pp_5)\cap S)\cup\{p_5\}$ and $S_3=Cone(p_3pp_5)\cap S$ or the partition $S_1=(R_2\cap S)\cup\{p, p_2, p_3\}$, $S_2=Cone(p_2pp_5)\cap S$ and $S_3=(Cone(p_3pp_5)\cap S)\cup\{p_5\}$, is admissible for $S$.

\item[{\it Case} 2.2] $p\notin R$. W.l.o.g. let $p\in\bbI(p_1p_2p_5)$. If $|R_1|+|R_5|\leq  3$, then $S_2=((R_1\cup R_5)\cap S)\cup\{p_1\}$ spans a $|S_2|$-hole or an empty $|S_2|$-pseudo-triangle. Moreover, $S_1=((R_2\cup R_4)\backslash \overline{\bbH}(p_3p_4, p_1))\cap S)\cup \{p, p_2, p_3, p_4, p_5\}$ spans a $|S_1|$-hole. Therefore, the partition of $S$ given by $S_1$, $S_2$, and $S_3=\overline{\bbH}(p_3p_4, p_1)\cap S$ is admissible. Otherwise, $|R_1|+|R_5|\geq 4$, which implies that $|R_1|=|R_5|=2$, from Observation \ref{ob:ob1}. If $|R_2|=0$, then $S_2=(R_5\cap S)\cup\{p, p_1, p_5\}$ forms a 5-hole and an admissible partition of $S$ is given by $S_1=(Cone(p_1pp_3)\cap S)\cup\{p_3\}$, $S_2$, and $S_3=Cone(p_5pp_3)\cap S$. Therefore, assume that $|R_2|\geq 1$. Then there exists a point $\beta \notin S$ on the line segment $p_2p_3$ such that $|Cone(p_1p\beta)\cap S|=|Cone(p_5p\beta)\cap S|=4$ (see Figure \ref{fig8}(c)), and the partition $S_1=(R_5\cap S) \cap\{p,p_1,p_5\}$, $S_2=Cone(p_1p\beta)\cap S$,
   and $S_3=Cone(p_5p\beta)\cap S$ is admissible.
\end{description}

\item[{\it Case} 3]$|L\{2, S\}|=6$. Consider the subdivision of the exterior of the second convex layer into 12
regions $R_i$ as shown in Figure \ref{fig7}(a). Note that the regions $R_i$ and $R_{i+2}$, for $i\in\{2, 4, 6, 8, 10, 12\}$ might intersect. Observation \ref{ob:ob1} implies that $S$ is admissible unless
$|\bigcup_{b=0}^2 R_{i+b}|\leq 2$, for $i\in\{2, 4, 6, 8, 10, 12\}$. Adding these inequalities and using the fact that $|\bigcup_{b=1}^{12} R_b|=7$, we get $|R_1|+|R_3|+|R_5|+|R_7|+|R_9|+|R_{11}|\geq 2$. This implies that for some $i\in \{1, 3, 5, 7, 9, 11\}$, $|R_i|\ne 0$. W.l.o.g. assume that $|R_7|\ne 0$. Let $Z_1=\mathcal H(p_4p_6, p_5)\cap S$ and $Z_2=\mathcal H(p_1p_3, p_2)\cap S$. From Observation \ref{ob:ob1} we know $|Z_1|\leq 5$ and $|Z_2|\leq 5$.

\item[{\it Case} 3.1] $|R_7|= 2$. If both $|Z_1|,|Z_2|\leq 4$, then $S_2=((R_1\cup R_7)\cap S)\cup\{p_1, p_3, p_4, p_6\}$ forms a 6-hole, and the partition $S_1=Z_1$, $S_2$, and $S_3=Z_3$ is an admissible partition of $S$. Next, assume that $|Z_1|=5$ and $|Z_2|= 2$. As $|R_7|=2$, $|R_6|=|R_8|=0$ by Observation \ref{ob:ob1}. Moreover, as $|Z_1|=5$, $|R_9|\ne 0$. If $|R_9|=2$, then
    the partition $S_1=((R_3 \cup R_9)\cap S)\cup\{p_1, p_2, p_4, p_5\}$, $S_2=\mathcal H(p_2p_4, p_3)\cap S$,
    and $S_3=\mathcal H(p_1p_5, p_6)\cap S$ is admissible. Otherwise, $|R_9|=1$ and $|\overline{\bbH}(p_5p_6, p_1)\cap S|=3$. The admissibility of $S$ then follows from Observation \ref{ob:ob1}.

\item[{\it Case} 3.2] $|R_7|=1$. In fact, by symmetry and {\it Case} 3.1 we can assume that $|R_i|\leq 1$ for all $i\in\{1, 3, 5, 9, 11\}$. If $|Z_1|=|Z_2|=4$, the partition $S_1=Z_1, S_2=((R_1\cup R_7)\cap S)\cup\{p_1, p_3, p_4, p_6\}, S_3=Z_2$ is admissible for $S$. Therefore, assume that $|Z_1|=5$ and $|Z_2|=3$. If $|R_8|+|R_9| \leq 1$, then $|\overline{\bbH}(p_5p_6, p_1)\cap S|\geq 3$ and admissibility follows from Observation \ref{ob:ob1}. As $|R_9|\leq 1$ by assumption, this implies that $|R_8|=|R_9|=1$. Therefore, $|R_{10}|=0$.
\begin{description}
\item[{\it Case} 3.2.1] $|R_{11}|=|R_{12}|=1$. In this case the partition of $S$ given by $S_1=((R_{11}\cup R_5)\cap S)\cup\{p_2, p_3, p_5, p_6\}$, $S_2=\mathcal H(p_2p_6, p_1)\cap S$, and $S_3=\mathcal H(p_3p_5, p_4)\cap S$ is admissible.  
\item[{\it Case} 3.2.2] $|R_{11}|=0$ and $|R_{12}|=2$. Then $|R_1|=|R_2|=0$.
\begin{description}
\item[{\it Case} 3.2.2.1] $|R_3|=1$. Then the admissible partition of $S$ is given by $S_1=((R_{9}\cup R_3)\cap S)\cup\{p_1, p_2, p_4, p_5\}$, $S_2=\mathcal H(p_1p_5, p_6)\cap S$, and $S_3=\mathcal H(p_2p_4, p_3)\cap S$.
\item[{\it Case} 3.2.2.2] $|R_3|=0$ but $|R_5|=1$. Then $S_1=((R_{11}\cup R_5)\cap S)\cup\{p_2, p_3, p_5, p_6\}$, $S_2=\mathcal H(p_2p_6, p_1)\cap S$, and $S_3=\mathcal H(p_3p_5, p_4)\cap S$ is admissible. 
\item[{\it Case} 3.2.2.3] $|R_3|=|R_5|=0$ and $|R_4|=2$. Observation \ref{ob:ob1} then implies that $|R_5|=|R_6|=0$ and  $|\overline{\bbH}(p_2p_3, p_6)\cap R_7\cap S|=0$. Therefore, $S_2=R_7\cup R_8 \cup \{p_2, p_3, p_4\}$ forms a 5-hole, which together $S_1=R_4\cup\{p_1, p_6\}$, and $S_3=R_9\cup R_{12}\cup\{p_5\}$ forms an admissible partition. 
\end{description}
\end{description}
\end{description}
\end{proof}

\subsection{$|CH(S)|=8$}

In this section we proof the admissibility of $S$ when $|CH(S)|=8$. As before assume that $L\{2, S\}=\{p_1, p_2, \ldots, p_m\}$, where the vertices are taken in counter-clockwise order, and the indices are to be identified modulo $m$.

\begin{lem}$S$ is admissible whenever $|CH(S)|=8$.
\label{lm:lmch8}
\end{lem}

\begin{proof}If $|L\{2, S\}|=3$ and none of the outer halfplanes induced
by the three edges of the second convex layer contains more than
two points of $\bbV(CH(S))$, then $|\bbV(CH(S))|\leq 2\times 3 < 8$.
Therefore, it suffices to assume that $4 \leq |L\{2, S\}| \leq 5$.

\begin{figure*}[h]
\centering
\begin{minipage}[c]{0.33\textwidth}
\centering
\includegraphics[width=1.95in]
    {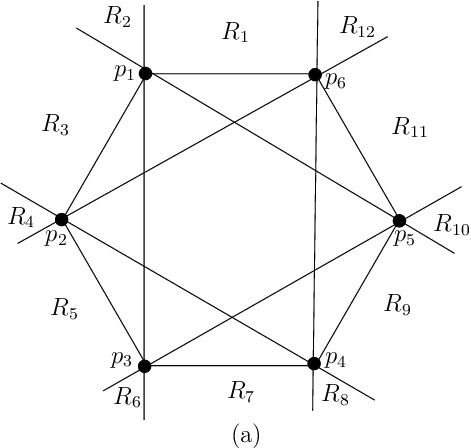}\\
\end{minipage}%
\begin{minipage}[c]{0.33\textwidth}
\centering
\includegraphics[width=1.80in]
    {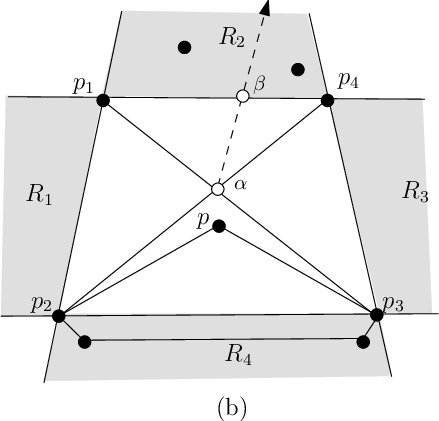}\\
\end{minipage}%
\begin{minipage}[c]{0.33\textwidth}
\centering
\includegraphics[width=1.95in]
    {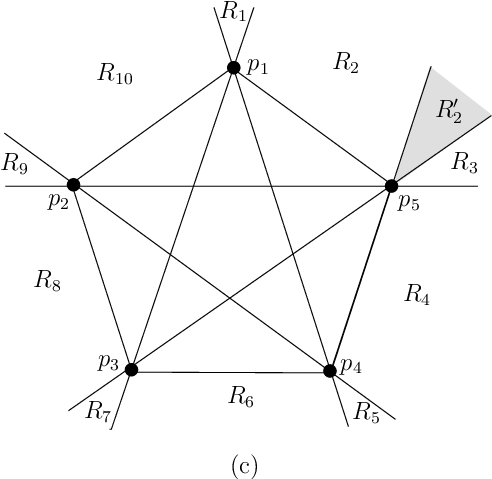}\\
\end{minipage}
\caption{Illustration for the proof of Lemma \ref{lm:lmch7} with $|L\{2, S\}| = 6$, and Illustration for the proof of Lemma \ref{lm:lmch8}: (b) $|L\{2, S\}| = 4$, (c) $|L\{2, S\}| = 5$.} 
\label{fig7}
\end{figure*}

\begin{description}
\item[{\it Case} 1]$|L\{2, S\}|=4$. Let $L\{3, S\}=\{p\}$, and $R_1, R_2, R_3, R_4$ be the 4 shaded regions outside the second convex layer as shown in Figure \ref{fig7}(b). Note that if $|R_1|+|R_3|\leq 3$, then $S_2=((R_1\cup R_3)\cap S)\cup\{p\}$ forms a convex quadrilateral or a 4-pseudo-triangle. Hence, $S_1=\overline{\bbH}_c(p_1p_4, p_2)\cap S$, $S_2$, and $S_3=\overline{\bbH}_c(p_2p_3, p_1)\cap S$ is an admissible partition of $S$. Therefore, we may assume that $|R_1|+|R_3|\geq 4$ and $|R_2|+|R_4|\geq 4$. This implies that $|R_{1}|=|R_{2}|=|R_{3}|=|R_{4}|=2$, as $\sum_{i=1}^{4}|R_i|\leq 8$. Let $\alpha$ be the point of intersection of the diagonals of the quadrilateral $p_1p_2p_3p_4$. W.l.o.g., assume that $p\in \bbI(p_2p_3\alpha)$. Then the points $p, p_2, p_3$ and along with two points in $R_{4}\cap S$ form a 5-hole. The remaining 8 points can be partitioned into two disjoint convex regions with 4 points each, because there exists a point $\beta$ on the line segment $p_1p_4$ such that $|Cone(\beta\alpha p_2)\cap S|=4$ and $|Cone(\beta\alpha p_3)\cap S|=4$ (see Figure \ref{fig7}(b)).

\item[{\it Case} 2]$|L\{2, S\}|=5$. Consider the partition of the exterior of the second convex layer into regions $R_i$'s as shown in Figure \ref{fig7}(c). Observation \ref{ob:ob1} implies that $S$ is admissible unless $|R_{i}|+|R_{i+1}|+|R_{i+2}|\leq 2$, for $i\in\{1, 3, 5, 7, 9\}$. Adding these 5 inequalities and using the fact that $\sum_{i=1}^{10}|R_i|=8$ we get, $|R_1|+|R_3|+|R_5|+|R_7|+|R_9|\leq 2$, that is, $|R_2|+|R_4|+|R_6|+|R_8|+|R_{10}|\geq 6$. Therefore, one of these 5 regions must contain exactly two points of $\bbV(CH(S))$. W.l.o.g., assume that $|R_{2}|=2$. Let $Z_1=\{\mathcal H(p_3p_5, p_4)\backslash R_7\}\cap \bbV(CH(S))$ and $Z_2=\{\mathcal H(p_3p_1, p_2)\backslash R_7\}\cap \bbV(CH(S))$.

\begin{description}
\item[{\it Case} 2.1] $|R_7|\geq 1$. We have $|R_5|+|R_6|\leq 1$ and $|R_8|+|R_9|\leq 1$. Therefore,
$|Z_1\cup\{p_4\}|=a_1\leq 4$ and $|Z_2\cup\{p_2\}|=a_2\leq 4$. Now, $|R_7|=6-(a_1+a_2)$ and there exists
a point $\alpha\in R_7\backslash S$ such that both $Cone(\alpha p_3p_1)\cap S$ and $Cone(\alpha p_3p_5)\cap S$ contain 4 points and $\{p_1, p_3, p_5\}\cup(R_2\cap\bbV(CH(S)))$ spans a 5-hole. Therefore, $S$ is admissible.

\item[{\it Case} 2.2] $|R_7|=0$. Let $R_{2}'=\overline{\bbH}(p_4p_5, p_1)\cap R_{2}$. We know that
$|Z_1|\leq 4$, $|Z_2|\leq 4$, and $|Z_1|+|Z_2|=6$. If $|Z_1|=|Z_2|=3$, the partition $S_1=Z_1\cup\{p_4\}$, $S_2=Z_2\cup\{p_2\}$, and $S_3=(R_2\cap S)\cup \{p_1,p_3,p_5\}$ is admissible. Otherwise, either $|Z_1|=4$ or $|Z_2|=4$. W.l.o.g. let $|Z_1|=4$. This implies that $|R_6|=2$, $|R_3|+|R_4|=2$ and $|R_2'|=0$. Let $s_i\in S$ be the first angular neighbor of $\overrightarrow{p_1p_3}$ in $R_6$. Then set $\{p_1,p_5, s_i\}\cup (R_2 \cap S)$ spans a 5-hole and the admissibility of $S$ follows. 
\end{description}
\end{description}
\end{proof}

\subsection{$|CH(S)|=9$}

The proof of Theorem \ref{lm3} can now be completed with the following simple observation:

\begin{obs}$S$ is admissible whenever $|CH(S)|\geq 9$.
\label{ob:obchgeq9}
\end{obs}

\begin{proof}Let $|CH(S)|=k \geq 9$. This implies that $|L\{2, S\}|=13-k\leq 4$, and so, there must exist an outer halfplane induced by an edge of $L\{2, S\}$ containing more than two points of $S$. The result now follows from Observation \ref{ob:ob1}. 
\end{proof}

\section{Conclusions}
\label{conclusion}
In this paper we prove that every set of 13 points, in general position, can be partitioned into three disjoint regions each of which span an empty convex polygon or an empty pseudo-triangle. This proves that the pseudo-convex partition number $\psi(n)\leq \lceil\frac{3n}{13}\rceil$, thus answering a question posed by Aichholzer et al. \cite{toth}.\\

\small{\noindent{\it Acknowledgement}: The authors are indebted to Bettina Speckmann for her insightful comments on the various properties of pseudo-triangles. The authors also thank the anonymous referees for valuable comments which improved the quality and the presentation of the paper.}

\end{document}